\documentclass[11pt,a4]{article}
\usepackage{amsmath,epsfig,fancyhdr, amsthm, amssymb, enumerate, mathtools}
\usepackage{fullpage}
\usepackage{graphicx}
\usepackage{subfigure}

\begin{document}\sloppy

\def\x{{\mathbf x}}
\def\L{{\cal L}}
\def\sampta{SampTA~}

\newcommand{\wt}[1]{\widetilde{#1}}
\newcommand{\mc}[1]{\mathcal{#1}}
\newcommand{\bigO}{\mathcal{O}}
\newcommand{\A}{\mathcal{A}}
\newcommand{\R}{\mathbb{R}}
\newcommand{\C}{\mathbb{C}}
\newcommand{\Z}{\mathbb{Z}}
\newcommand{\N}{\mathbb{N}}
\newcommand{\E}{\mathbb{E}}
\renewcommand{\P}{\mathbb{P}}
\newcommand{\sA}{\mathscr{A}}
\newcommand{\sM}{\mathscr{M}}
\newcommand{\sY}{\mathscr{Y}}
\newcommand{\sD}{\mathscr{D}}
\newcommand{\sG}{\mathscr{G}}
\newcommand{\sd}{\Sigma\Delta}
\newcommand{\psum}{\mathop{\sum\nolimits'}}
\newcommand{\pprod}{\mathop{\prod\nolimits'}}
\newcommand{\diag}{\operatorname{diag}}

%
%

\newtheorem{theorem}{Theorem}[]
\newtheorem{lemma}[theorem]{Lemma}
\newtheorem{proposition}[theorem]{Proposition}
\newtheorem{cor}[theorem]{Corollary}
\theoremstyle{remark}
\newtheorem{remark}{Remark}[theorem]
\theoremstyle{definition}
\newtheorem{definition}[theorem]{Definition}

\title{Sparse recovery for spherical harmonic expansions}

\author{Holger Rauhut\thanks{Hausdorff Center for Mathematics \& Institute for Numerical Simulation, Universit{\"a}t Bonn} and Rachel Ward\thanks{Courant Institute, NYU, 251 Mercer Street, NY, NY, 10003, U.S.A.}}

\maketitle

\begin{abstract}
We show that sparse spherical harmonic expansions can be efficiently recovered from a small
number of randomly chosen samples on the sphere.  To establish the main result, we verify the restricted isometry property of an associated preconditioned random measurement matrix using recent estimates on the uniform growth of Jacobi polynomials.
\end{abstract}

\section{Introduction}

Compressed sensing has triggered significant research activity in recent years. It predicts that sparse signals can be recovered from what was previously believed to be highly incomplete information. In this work, we show that functions on 
the sphere $S^2 = \{\x \in \R^3, \|\x\|_2 = 1\}$ that have a sparse or compressible representation in 
the spherical harmonic basis can be recovered from a number of samples that scales (essentially) 
linearly with the sparsity level.  This  can be viewed as an extension of existing results \cite{2,6} 
for sparse recovery of trigonometric polynomials on the circle.  
Since the ($L_2$-normalized) spherical harmonic basis functions are not uniformly bounded, standard compressed sensing 
theory does not apply directly to sparse recovery on the sphere.  Instead, we appeal to recent results in \cite{1} concerning sparse recovery in orthonormal polynomial systems.  Since orthonormal polynomials blow up sufficiently quickly and uniformly at the endpoints of their domain, one preconditions in order to transform the problem into that of 
sparse recovery in a uniformly bounded system.  The decomposition of spherical harmonic basis functions
into tensor products of trigonometric polynomials and Jacobi polynomials 
allows to prove the main result: any degree-$D$ polynomial on the sphere (that is, with $N = D^2$ coefficients) 
consisting of at most $s$ spherical harmonic basis elements can be efficiently recovered from $m \sim s N^{1/4} \log^4{(N)}$ independent sampling points drawn uniformly with respect to a certain measure (see below).
We establish this result by verifying the restricted isometry property (RIP) of an associated random matrix. 

\section{Background and notation}

The spherical harmonics $Y_{\ell}^k, -\ell \leq k \leq \ell, k \geq 0$ form an orthonormal basis for the Hilbert space of square-integrable functions on the sphere.  They are orthogonal with respect to the spherical surface measure $\Omega$. In spherical coordinates 
$(\phi, \theta) \in B := [0,\pi] \times [0,2\pi)$, $(x = \cos(\theta)\sin(\phi), y = \sin(\theta)\sin(\phi), z = \cos(\phi)) \in S^2$, this orthogonality relation becomes
\begin{equation}
\label{eq:orthonormal}
\int_{0}^{2\pi} \int_{0}^{\pi} 
Y_{\ell}^k (\phi, \theta) \overline{Y_{\ell'}^{k'} (\phi, \theta)} \sin(\phi) d\phi d\theta = \delta_{\ell \ell'} \delta_{k k'}.
\end{equation}
Here, $\overline{z}$ denotes the complex conjugate and $\delta_{mn}$ is the Kronecker delta.
The spherical harmonics may be expressed as 
\begin{eqnarray}
\label{eq:harmonics}
Y_{\ell}^k (\phi, \theta) = e^{ik\theta} (\sin{\phi})^{|k|} p_{\ell - |k|}^{|k|} (\cos{\phi}),\;\; (\phi, \theta) \in B,	
\end{eqnarray}
where the Jacobi polynomials $(p_{n}^{\alpha})_{n=0}^{\infty}$ with parameter $\alpha \geq -1$ are 
the orthonormal polynomial basis on the interval $[-1,1]$ with respect to the measure $dv(x) = (1 - x^2)^{\alpha} dx$  \cite{4}. 
In particular,  the Lebesgue measure $\alpha = 0$ generates the Legendre polynomials, while the Chebyshev 
measure, $\alpha = -1/2$, generates the Chebyshev polynomials.  

Spherical harmonic expansions on the sphere are analogous to Fourier series expansions on the circle.  
Functions on the sphere of the form 
\begin{equation}
\label{eq:trig}
g(\phi, \theta) = \sum_{\ell=0}^{D-1} \sum_{k=-\ell}^{\ell} c_{\ell, k} Y_{\ell}^k (\phi, \theta)
\end{equation}
are called \emph{harmonic polynomials} of degree $D-1$.  Note that $N = D^2$ spherical harmonic basis elements generate harmonic polynomials of degree $D-1$.  We will call a harmonic polynomial \emph{s-sparse} if its coefficient vector $c = (c_{\ell, k}) \in \mathbb{C}^N$ has cardinality at most 
$s$; i.e. $\| c \|_0 := | \{(\ell, k): c_{\ell, k} \neq 0 \} | \leq s$. More generally, the degree to which a harmonic 
polynomial can be well-approximated by its $s$ most significant coefficients can be quantified using the concept of best
 $s$-term approximation error 
 which is defined, 
for a vector $z \in \mathbb{C}^N$ by 
$$
\sigma_s(z)_1 = \inf_{y: \| y \|_0 \leq s} \| y - z \|_1 .
$$
We say that a harmonic polynomial \eqref{eq:trig}  is \emph{compressible} if $\sigma_s(c)_1$ decays quickly as $s$ increases.

\section{Main results}

We aim to recover sparse harmonic polynomials on the sphere from only a few function samples.
Note that samples $y_j = g(\phi_j, \theta_j)$,  $j=1,\hdots,m$, of a $D-1$-degree harmonic polynomial 
$g(\phi, \theta) = \sum_{\ell=0}^{D-1} \sum_{k=-\ell}^{\ell} c_{\ell, k} Y_{\ell}^k (\phi, \theta)$ may be expressed concisely in terms of the coefficient vector $c = (c_{\ell, k}) \in \mathbb{C}^{N}$ according to
\begin{equation}
\label{eq:linear1}
y = 
\Phi c.
\end{equation}
where $\Phi$ is the $m \times N$ matrix defined component-wise by 
\begin{equation}
\label{eq:shmatrix}
\Phi_{j,(k,\ell)} =Y_{\ell}^k (\phi_j,\theta_j).
\end{equation}
We are interested in solving the system of linear equations \eqref{eq:linear1}  in the underdetermined setting $m < N$, and in particular, to single out the original sparse coefficient vector $c$ from among the infinitely-many solutions. The compressed sensing literature has suggested various  reconstruction algorithms for sparse recovery; for simplicity we focus only on $\ell_1$-minimization \cite{7} in this paper.  

The $L_{\infty}$-norm of the  spherical harmonics $Y_{\ell}^k(\phi, \theta)$ increases with the degree 
$\ell$ according to $\sup_{-\ell \leq k \leq \ell} \| Y_{\ell}^k \|_{\infty} = \sup \| Y_{\ell}^0 \|_{\infty} = \ell$, and this extremum is obtained at the spherical caps $\phi = 0, \pi$.  
This means that the linear system \eqref{eq:linear1} is in general ill-conditioned.   
It is known (see Proposition \ref{eq:spherebound}) that $| (\sin{\phi})^{1/2} Y_{\ell}^k(\phi, \theta) | \sim \ell^{1/4}$; consequently, we precondition the system \eqref{eq:linear1} for numerical stability, multiplying both sides by  the $m \times m$ diagonal matrix ${\cal{A}}$ with entries ${\cal{A}}_{j,j} = (\sin{\phi_j})^{1/2}$,
\begin{equation}
\label{eq:linear2}
{\cal A} y = 
{\cal A} \Phi c.
\end{equation}

Our main result is that any $s$-sparse harmonic polynomial on the sphere of maximal degree $D$ can be recovered efficiently from a number of samples $m$ that scales linearly with the sparsity level and sublinearly with the degree.  This reconstruction is moreover robust with respect to noisy samples and passing from sparse to compressible vectors.

\begin{theorem}
Let $m, s$, and $N$ be given integers satisfying
 \begin{equation}
 \label{eq:m}
m \geq C s \log^3(s) N^{1/4}\log(N).
 \end{equation}
Suppose that $m$ coordinates on the sphere $(\phi_1, \theta_1), \hdots, (\phi_m, \theta_m)$ are drawn independently from the uniform measure on $B = [0,\pi] \times [0,2\pi)$. 

Let $\Phi$ be the $m \times N$ spherical harmonic matrix \eqref{eq:shmatrix} and let $ {\cal{A}} \Phi$ 
be its preconditioned version \eqref{eq:linear2}.  With probability exceeding $1-N^{-\gamma \log^3(s)}$ the following holds for all harmonic polynomials $g(\phi, \theta) = \sum_{\ell=0}^{D-1} \sum_{k=-\ell}^{\ell} c_{\ell, k} Y_{\ell}^k(\phi, \theta)$.  Suppose that noisy sample values $y_j = g(\phi_j, \theta_j)  + \eta_j$ are observed, and that $\| \eta \|_{\infty} \leq \varepsilon$.  
Let
\begin{equation}
\label{ell1}
\widehat{c} = \arg \min \| z \|_1 \hspace{3mm} \textrm{subject to} \hspace{3mm} \| {\cal A}\Phi z - {\cal A} y \|_2 \leq \sqrt{m}\varepsilon.
\end{equation}
Then
$$
\| c - \widehat{c} \|_2 \leq \frac{C_1 \sigma_s(c)_1}{\sqrt{s}} + C_2\varepsilon.
$$
The constants $C, C_1, C_2,$ and $\gamma$ are universal.
\end{theorem}

We expect that the bound \eqref{eq:m} is not optimal, but so far we have not been able to remove the polynomial factor $N^{1/4}$.

\section{Numerical experiments}

In Figure $1$, we plot two phase diagrams illustrating the success of $\ell_1$-minimization (i.e., Equation \eqref{ell1} with $\varepsilon = 0$) in recovering $s$-sparse harmonic polynomials $g(\phi, \theta) = \sum_{\ell=0}^{D-1} \sum_{k=-\ell}^{\ell} c_{\ell, k} Y_{\ell}^k (\phi, \theta)$ from $m$ sample values $g(\phi_j, \theta_j)$.    In each plot, we fix the degree to be $D = N^{1/2} = 16$ and vary the sparsity level $s$ and number of measurements $m$.  We form an $s$-sparse coefficient vector $c = (c_{\ell, k} )$ by choosing a random support set from $\mathbb{N} \cap [1, N]$ of cardinality $s$ and assigning independent and identically distributed Gaussian weights as the coefficients to this support;  we then draw $m$ sampling points
which we use to recover the $s$-sparse vector using $\ell_1$-minimization.  For each pair $(s,m)$, we record the frequency of success of $\ell_1$-minimization out of $20$ trials.   In Figure $1(a)$ the sampling points are chosen independently from the product measure $d\phi d\theta$ on $B$, which has higher sampling density around the spherical poles and for which the sparse recovery results of Theorem $1$ apply.   We observe a large region of phase space (in black) corresponding to uniform recovery, similar to that observed in the phase transition curves obtained for other compressive sensing matrices \cite{5}.  On the other hand, in Figure $1(b)$ the sampling points are chosen independently from the uniform surface measure $\sin{(\phi)} d\phi d\theta$.  In this case, $\ell_1$-minimization fails to recover sparse harmonic polynomials for essentially all parameters $(s,m)$.

\begin{figure}
\label{fig:1}
\begin{center}
\caption{Phase transition for sparse recovery on the sphere}
\subfigure[Sampling points from product measure $d\phi d\theta$]{\label{sugfig1}\includegraphics[width=7cm]{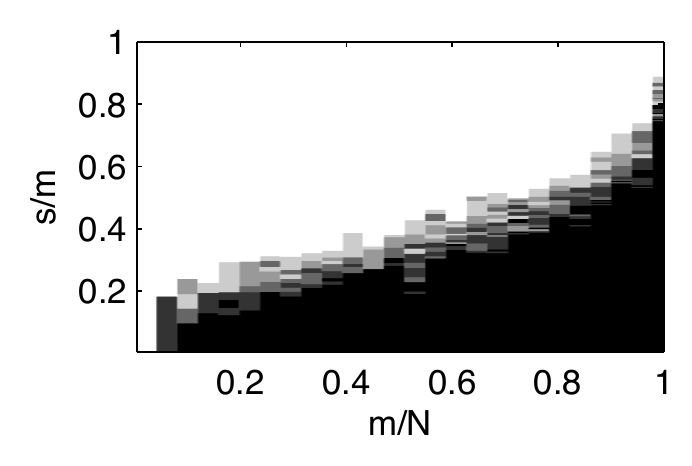} }
\subfigure[Uniformly distributed sampling points]{\label{subfig2} \includegraphics[width=7cm]{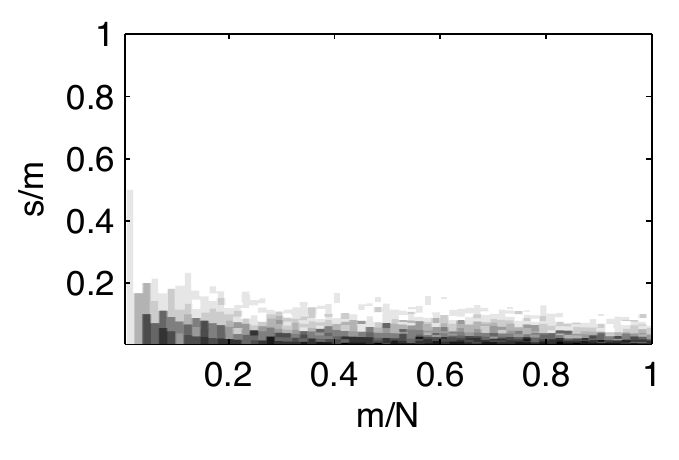}}

\end{center}
\end{figure}

\section{Sparse recovery via Restricted Isometry Constants}

 We prove Theorem $1$  by showing that the preconditioned spherical harmonic matrix $\widetilde{\Phi} = {\cal A} \Phi$ in \eqref{eq:linear2} satisfies the  \emph{restricted isometry property} (RIP) \cite{2,7}.  

\begin{definition}[Restricted isometry constants]
Let $\Psi \in \mathbb{C}^{m \times N}$.  For $s \leq N$, the restricted isometry constant $\delta_s$
associated to $\Psi$ is the smallest number such that
\begin{equation}\label{def:RIP}
(1-\delta_s) \|c\|_2^2 \leq \|\Psi c\|_2^2 \leq (1+\delta_s) \|c\|_2^2
\end{equation}
for all $s$-sparse vectors $c \in \C^N$.
\end{definition}

Informally, we say that the matrix $\Psi$ ``has the restricted isometry property" if
$\delta_s$ is small for $s$ reasonably large compared to $m$.
For matrices satisfying the restricted isometry property, the following $\ell_1$-recovery results can be shown
\cite{7,8}.

\begin{theorem}[Sparse recovery for RIP-matrices]
\label{thm:l1:stable} 
Let $\Psi \in \C^{m \times N}$. Assume
that the restricted isometry constant of $\Psi \in \C^{m \times N}$
satisfies
\begin{equation}
\label{RIP:const}
\delta_{2s} < 3/(4 + \sqrt{6}) \approx 0.4652.
\end{equation}
Let $x \in \C^N$ and assume noisy measurements $y = \Psi x + \eta$ are given with $\|\eta\|_2 \leq \varepsilon$. Let $x^\#$ be the minimizer of 
\begin{align}\label{l1eps:prog}
\arg \min_{z \in \C^N} \quad \| z \|_1 \mbox{ subject to } \quad \|\Phi z - y \|_2 \leq \varepsilon.
\end{align}
Then
\begin{align}
\label{l2noise}
\|x - x^\#\|_2 \leq  C_1 \frac{\sigma_s(x)_1}{\sqrt{s}} + C_2 \varepsilon,
\end{align}
for some constants $C_1, C_2 >0$ that depend only on $\delta_{2s}$.
In particular, if $x$ is $s$-sparse then reconstruction is exact, $x^\# = x$.
\end{theorem}

A general setup for matrices having the restricted isometry property are those 
associated to \emph{bounded orthonormal systems} \cite{2,6,7}.


\begin{theorem}[RIP for bounded orthonormal systems]
\label{thm:BOS:RIP} 
Consider an orthonormal system of functions $\psi_j$, $j \in [1,N] \cap \N$ on a measurable space ${\cal M}$ endowed with a probability measure $\nu$,
that is $\int_{\cal M} \psi_j \overline{\psi_k} d\nu = \delta_{j,k}$. 
Consider the matrix $\Psi \in \C^{m \times N}$ with entries
\begin{equation}\label{def:Phi:matrix}
\Psi_{\ell,k} = \psi_k(x_\ell), \quad \ell \in [1,m] \cap \N, k \in [1,N] \cap \N,
\nonumber
\end{equation}
formed by i.i.d.\ samples $x_\ell$ drawn from the measure $\nu$.  Suppose this system 
has the uniform bound $K = \sup_{j \in [N]} \|\psi_j\|_\infty = \sup_{j \in [N]} \sup_{x \in {\cal D}} |\psi_j(x)|$.
If
\begin{equation}\label{BOS:RIP:cond}
m \geq C\delta^{-2} K^2 s \log^3(s) \log(N),
\end{equation}
then with probability at least 
$1-N^{-\gamma \log^3(s)},$ 
the restricted isometry constant $\delta_s$ of $\frac{1}{\sqrt{m}} \Psi$ 
satisfies $\delta_s \leq \delta$. The constants $C > 0$ and $\gamma>0$ are universal.
\end{theorem}

An important special case  is the matrix associated to samples of the trigonometric system $(\exp(2 \pi i j x))_{j=0}^{N-1}$ chosen from the uniform measure on $[0,1]$, which has the optimal uniform bound $K=1$.  
Another example is the sampling matrix associated to the Chebyshev polynomial system. 
In this case, $K = \sqrt{2}$.  

\section{Sparse recovery in spherical harmonic systems}

Recall from \eqref{eq:harmonics} that the spherical harmonics can be expressed as tensor products of complex exponentials in $\theta$ and orthogonal polynomials in $\cos(\phi)$.  Since the latter are not uniformly bounded, spherical harmonics do not fall directly into the scope of bounded orthonormal systems.   To get around this obstacle, we proceed in a similar fashion to \cite{1}, and use estimates on the uniform rate of growth of orthogonal polynomials in order to precondition 
the spherical harmonic system. 

First we will need the following growth estimates.
\begin{proposition}
\label{thm:kras}
Consider the weight function $v(x) = (1 - x^2)^{\alpha}$ on $[-1, 1]$, and let $(p^{\alpha}_n)_{n}$ be the associated orthonormal polynomial system.  Then, for all $x \in [-1,1]$, the following holds.
\begin{enumerate}\itemsep-2pt
\item If $\alpha = 0$, the associated polynomials $(p_n^0)$ are the Legendre polynomials and satisfy
\begin{equation}
\label{legendre}
(1 - x^2)^{1/4} |p_n^{0}(x) | \leq 2/\sqrt{\pi}.
\end{equation}
\item For any $\alpha \geq 0$,
\begin{equation}
\label{bound:classical}
(1 - x^2)^{1/4 + \alpha/2} |p_n^{\alpha}(x) | \leq C_{\alpha}
\end{equation}
\item  If $\alpha \geq 3$, then
\begin{equation}
\label{eq:ultrabound}
(1 - x^2)^{1/4 + \alpha/2} |p_n^{\alpha}(x) | \leq C \alpha^{1/6} \Big( 1 + \frac{\alpha}{n} \Big)^{1/12},
\end{equation}
where $C > 0$ is a universal constant.
\end{enumerate}
\end{proposition}
The bound \eqref{legendre} for Legendre polynomials is classical and known to be tight, and the more general bound \eqref{bound:classical} is also classical; see \cite{4} for more details.    The more refined bound \eqref{eq:ultrabound} was derived only recently in \cite{3}.  Although the result is stated in \cite{3} for the parameter range $\alpha \geq 1, n \geq 6$, it is not hard to verify that the key estimate (Lemma $8$ in \cite{9}) is valid also for the parameter range $\alpha \geq 3, n \geq 0$.

Using these bounds in conjunction with the tensor product representation \eqref{eq:harmonics}, we arrive at the following rate of growth for the spherical harmonics.
\begin{proposition}
\label{eq:spherebound}
For all $\ell \in \Z^+$ and $-\ell \leq k \leq \ell$,
\begin{align}
(\sin{\phi})^{1/2} | Y_{\ell}^k (\phi, \theta) | \leq C (\ell+1)^{1/4} \hspace{3mm} \forall x \in [-1,1]. \nonumber
\end{align}
\end{proposition}

We can now state the proof of Theorem $1$.

\vspace{3mm}

\noindent { \bf Proof of Theorem $1$.}

\vspace{2mm} 
Consider, for  $0\leq \ell \leq D-1$, $-\ell \leq k \leq \ell$, the functions 
\begin{equation}
\label{Q}
Q_{\ell}^k (\phi, \theta) = (\sin{\phi})^{1/2} Y_{\ell}^k (\phi, \theta). 
\end{equation}
By Proposition \ref{eq:spherebound}, 
$$\sup_{0\leq \ell \leq \sqrt{N}-1 , -\ell \leq k \leq \ell}  \| Q_{\ell}^k \|_{\infty} \leq C N^{1/8} $$
 for a universal constant $C$.  Because the spherical harmonics $Y_{\ell}^{k}$ are orthonormal with respect to the uniform measure $\sin(\phi)d\phi d\theta$, the $Q_{\ell}^k$'s are orthonormal with respect to the product measure $d\theta d\phi$:
\begin{align}
\label{ortho_nu}
&\int_{0}^{2\pi} \int_{0}^{\pi} Q_{\ell}^{k}(\phi, \theta) Q^{k'}_{\ell'}(\phi, \theta) d\theta d\phi \\ 
&=  \int_{0}^{2\pi} \int_{0}^{\pi} Y_{\ell}^{k}(\phi, \theta) Y^{k'}_{\ell'}(\phi, \theta) \sin(\phi) d\theta d\phi  
= \delta_{\ell \ell'} \delta_{k k'}. \nonumber
\end{align}
Applying Theorem \ref{thm:BOS:RIP} to the system $\{ Q_{\ell}^k \}$ , whose sampling matrix $Q_{\ell}^k(\phi_j, \theta_j)= y_j$ is equivalent to the preconditioned spherical harmonic matrix $(\sin{\phi_j})^{1/2} Y_{\ell}^k (\phi_j, \theta_j) = y_j$ \eqref{eq:linear2},
Theorem $1$ follows from the recovery results for restricted isometry systems in Theorem \ref{thm:l1:stable}.

\bibliographystyle{plain}

\end{document}